\documentclass[10pt]{article}
\usepackage[utf8]{inputenc}
\usepackage[T1]{fontenc}
\usepackage[english]{babel}
\usepackage{amsmath, setspace, amsfonts, amssymb, graphics,multirow}
\usepackage{amsthm}
\usepackage{graphicx}
\usepackage{natbib}
\usepackage{color}
\usepackage{soul}
\usepackage{cancel}

\usepackage{authblk}
\usepackage[paperwidth=19cm, paperheight=28cm]{geometry}
\makeatletter

\newcommand{\diag}{\mathop{\mathrm{diag}}}

\newcommand{\bbeta}{\boldsymbol{\beta}}

\newcommand{\btheta}{\boldsymbol{\theta}}
\newcommand{\bpsi}{\boldsymbol{\psi}}

\newcommand{\bC}{\mathbf{C}}

\newcommand{\bZ}{\mathbf{Z}}
\newcommand{\bu}{\mathbf{u}}
\newcommand{\by}{\mathbf{y}}
\newcommand{\bz}{\mathbf{z}}
\renewcommand{\bf}{\mathbf{f}}
\newcommand{\bt}{\mathbf{t}}
\newcommand{\bm}{\mathbf{m}}
\newcommand{\br}{\mathbf{r}}

 \DeclareMathOperator{\Tr}{Tr}

\newtheorem{proposition}{Proposition}

\newtheorem{algo}{\textit{Algorithm}}

\author[1,2]{Pierre Barbillon}
\author[3]{C\'elia Barth\'el\'emy}
\author[4,5]{Adeline Samson}
\affil[1]{AgroParisTech / UMR INRA MIA, F-75005 Paris}
\affil[2]{INRA, UMR 518, F-75005 Paris}
\affil[3]{INRIA Saclay, Popix team, Orsay}
\affil[4]{Univ. Grenoble Alpes, LJK, F-38000 Grenoble, France}
\affil[4]{CNRS, LJK, F-38000 Grenoble, France}

\title{Parametric estimation of complex mixed models based on meta-model approach}
\begin{document}
 \maketitle

\begin{abstract}
Complex biological  processes are usually experimented  
along time among a collection of individuals. Longitudinal data are then available and the statistical challenge is  to better understand the underlying biological mechanisms. 
The standard   statistical approach  is mixed-effects model, with  regression functions that 
are now highly-developed to describe precisely the biological processes (solutions of
multi-dimensional ordinary differential equations or of partial differential equation).
When there is no analytical solution,
a classical estimation approach   relies on the coupling of a stochastic version of the EM algorithm (SAEM) 
with a MCMC algorithm. This procedure needs many evaluations of the regression function 
which is clearly prohibitive when a time-consuming solver is used for computing it.
In this work a meta-model relying on a Gaussian process emulator is proposed to replace this regression function.
The new source of uncertainty due to this approximation can be incorporated in the model which leads
to what is called a mixed meta-model.
A control on the distance 
between the maximum likelihood estimates in this mixed meta-model
and the maximum likelihood estimates obtained with the exact mixed model is guaranteed.
Eventually, numerical simulations are performed to illustrate the efficiency of this approach.\\
\smallskip

\end{abstract}

\section{Introduction}
\label{intro}

Mixed-effects model methodology \citep{Pinheiro2000}
  allows to discriminate between inter- and intra-subjects variabilities which is essential when dealing with longitudinal data. 
  Statistical methods for mixed models are now well established (see references below) but can be  time consuming depending on the complexity of the regression functions. Indeed sophisticated
mathematical models have been developed to describe precisely biological processes: multi-dimensional ordinary differential equations (ODE) \citep[see][for modeling viral load decrease in HIV patients or tumor growth]{Wu2005, Guedj2007, Lavielle2011, Ribba2012}  or   partial differential equation (PDE)  \citep[see][for modeling   tumor growth or HVC viral kinetic]{Grenier2014, Chatterjee2012}. 
 These mathematical models have no analytical solution and only an approximate solution can be obtained with computationally intensive numerical methods. 
  This   induces a huge increase of the  computation cost of the   estimation method \citep[up to 23 days according to][]{Grenier2014}. Therefore, there is a crucial need to develop  new statistical approaches to reduce the computation time. 
%
%\medskip
The significant computation time of mixed models estimation methods  is due to their iterative  settings, compulsory to sidestep the intractability of the likelihood. This is true for methods based on   linearisation  \citep{Pinheiro2000} or likelihood numerical approximation  \citep{Davidian1995, Wolfinger1993} and this is crucial for EM-type methods such as stochastic   EM algorithms  \citep{Wei1990, Kuhn2005}. 

\medskip
Our objective   is to propose a way of reducing the computation time of the SAEM-MCMC algorithm \citep{Kuhn2005} for complex mixed models, together with a theoretical study of the resulting estimator. 
When the regression function is not analytically available (and we call it also computer model in the rest of the paper), extensions of SAEM have already been proposed.  \cite{Donnet2007}   deal with the case of an ODE mixed model, approximating the solution with a numerical scheme and studying the influence of this scheme  
to the properties of the approximate maximum likelihood estimator. But this approach remains time consuming when the ODE is multi-dimensional. For a PDE mixed model,   \cite{Grenier2014}  propose to approximate the PDE  with a numerical scheme on a predefined grid, and then to interpolate  the solution of the PDE linearly between two points of the grid.
This linear approximation allows  substantially reducing the computation time from 23 days to around 30 minutes, but may lead to biased estimates depending on the non linearity of the model.
% Although these first results are very promising, the theoretical properties of the obtained estimator  has not been studied.
%, which will 
%\textbf{Voir aussi si on peut citer le premier travail de Celia ???}

\medskip
The keystone of providing an efficient estimation method with good statistical properties is the choice of the procedure approximating 
the regression function. 
More accurate surrogates of computer model than linear approximation rely on Gaussian process emulation which consists of modeling the computer model 
as the realization of a Gaussian process \citep{sacks:schiller:1989,santner:etal:2003,fang:2005}.
This technique is  also known as Kriging.
A cheap approximation of the model, the emulator, is obtained by conditioning the Gaussian process on some evaluations of the model
corresponding to inputs of a well-chosen design of numerical experiments.
This stochastic modeling of the computer model provides also a measure of uncertainty on the precision of the approximation 
as a supplementary variance-covariance function which can be integrated in the mixed model.
This approach has been already coupled with a Stochastic EM algorithm \citep{barbillon:2011} or with a Bayesian
procedure  \citep{fu:2014} in regression models without random effects.
In this paper, we propose to couple the SAEM algorithm with this Gaussian process emulator,     incorporating this new source of uncertainty due to the approximation. 
Thus, providing a confidence interval   of the unknown parameters   takes into account the error induced by the approximation.
We will refer to this approach as the complete mixed meta-model.
However, the supplementary variance-covariance function 
in the model induces a loss of independence of the observations obtained from the different subjects which
increases the computational burdensome of the MCMC scheme.
That is why we also propose two simplified versions: the first one (called intermediate) by considering a diagonal variance-covariance function of the approximation error; the second one (called simple)  by only using the approximation  of the computer model and not incorporating the
variance-covariance function.
These two last versions allow to assume the independence between the subjects, and to reduce significantly the computational cost. %This is illustrated convincingly on simulated examples. 

\medskip
The Gaussian process emulator can also be interpreted as an approximation of the computer model
by kernel interpolation with radial basis function as in \cite{schaback:1995,schaback:2007}.
%It consists in the same approximation as Gaussian process emulator.
% The computer model can also be approximated by  kernel interpolation with radial basis function as in \cite{schaback:1995,schaback:2007}, leading 
% to the same approximation as the Gaussian process emulation.
In this framework,
a point-wise control on the error of approximation is provided. Hence, we are able to guarantee a control on the distance 
between the maximum likelihood estimates in the approximate mixed meta-models
and the maximum likelihood estimates obtained with the exact computer model.
This control is decreasing to zero as a function of the space-fillingness of the design of numerical experiments.

 \medskip
The paper is organized as follows. 
Section \ref{sec:MixedModel} introduces the standard non-linear mixed model and Section
\ref{sec:MetaModel} recalls the principles and the main results of the Gaussian process emulation. 
Section \ref{sec:MixedMetaModel} introduces three mixed models approximated by Gaussian process emulator. 
In Section \ref{sec:estimation},  three versions of the SAEM algorithm coupled to a Gaussian process emulator are proposed. 
Theoretical results are given in Section \ref{sec:ConvApproxMLE}.  A simulation study illustrates these results 
(Section \ref{sec:simulation}). Section \ref{sec:conclusion} concludes the paper with some extensions. Proofs are gathered in Appendix.

\section{Mixed model and notations}\label{sec:MixedModel}
Let us define $\by_i=\,^t(y_{i1},\ldots,y_{in_i})$ where $y_{ij}\in\mathbb{R}^p$ is the response  for individual $i$ at time $t_{ij}$, $i=1,\ldots,N$, $j=1,\ldots,n_i$, and let $\by=(\by_1,\ldots,\by_N)$ be the vector of all observations, of size $n_{tot}=\sum_{i=1}^N n_i$. We assume that the individual vectors $\by_{i}$ are 
described by a non-linear mixed model, defined as follows, for $j=1,\ldots,n_i$:
 \begin{eqnarray}\label{eq:MixedModel}
y_{ij} &= &  f (t_{ij}, \psi_i) + \sigma_\varepsilon \, \varepsilon_{ij}, \quad  \varepsilon_{ij} \sim_{iid} \mathcal{N}(0,1)\,,\\
 \psi_i &\sim_{iid}&    \mathcal{N}(\mu,\Omega), \nonumber
 \end{eqnarray}
where  $f(\cdot, \cdot): \mathbb{R}\times \mathbb{R}^d\rightarrow\mathbb{R}^p$  is the regression  function,  %and  $G(\cdot):   \mathbb{R}^p\rightarrow\mathbb{R}^q$  is the function of the observation model (we allow $q<p$),    
$\psi_i$ is a $d$-vector of individual parameters. The $\varepsilon_i=(\varepsilon_{i1},\ldots,\varepsilon_{in_i})^t$ represents the Gaussian centered residual error, independent of $\psi_i$. The individual parameter $\psi_i$ are assumed to be random, Gaussian with expectation  $\mu$    and $d\times d$ covariance matrix $\Omega$. Note that the individual vectors $(\by_{i})_{i}$ are independent and identically distributed. 

 The quantities we want to estimate from the observations $\by$ are the population parameters $\theta = (\mu, \Omega, \sigma_\varepsilon^2)$.  In the following, we restrict to the case of scalar observations ($p=1$) to ease the reading, but the extension to a multidimensional observation with $p>1$ is straightforward. 

 \medskip

We want to estimate $\theta$ by maximum likelihood. The likelihood of model (\ref{eq:MixedModel}) is:\begin{eqnarray}\label{eq:Likelihood}
&&p(\by,\theta) = \int p(\by,\bpsi\, ; \theta)\,d\bpsi = \prod_{i=1}^N\int p(\by_i|\psi_i; \theta) p(\psi_i; \theta) d\psi_i
\nonumber\\
&=&\prod_{i=1}^N \int 
\Bigg\{\frac{1}{(2\pi \sigma_\varepsilon^2)^{n_i/2}} 
\times
\exp\left(-\frac{1}{2}\,^t(\by_i-\bf(\bt_i,\psi_i))(\sigma_\varepsilon^2 I_{n_i})^{-1}(\by_i-\bf(\bt_i,\psi_i)) \right)\nonumber \\
&&\times
\frac1{(2\pi)^{d/2}|\Omega|^{1/2 }} \exp\left( -\frac{1}{2}  \,^t(\psi_i -\mu) \Omega^{-1} (\psi_i - \mu)\right) d\psi_i\Bigg\}
\end{eqnarray}
where $\bt_i=(t_{i1}, \ldots, t_{in_i})$ and\\ $\bf(\bt_i,\psi_i)=\,^t(f(t_{i1}, \psi_i), \ldots, f(t_{in_i}, \psi_i))$. 
When $f$ is non linear with respect to $\bpsi$, the maximum likelihood estimator has no closed form. %Thus, there is a need to sidestep this first difficulty. But 
Any estimation method adapted to non-linear mixed models would require a very large number of evaluations of $f$, which could be time consuming when the structural function  $f$ is a computer model. 
Therefore, there is a real need to consider approximations of the function $f$ that are simple to evaluate at any point $(t, \psi)$. For that purpose, we introduce the framework of meta-model in  the next section.

%%%%%%%%%%%%%  
\section{Meta-model}\label{sec:MetaModel} 
 We start with the point of view of conditioned Gaussian process
 emulation which has the nice feature of incorporating as a variance-covariance function
 the additional source of uncertainty due to the approximation.
 This will naturally lead to a mixed meta-model 
 on which the SAEM-MCMC can be performed.
 We also link this framework to the kernel interpolation framework
 since we need the deterministic point-wise control on the approximation to obtain the control between
 the maximum likelihood estimates corresponding to the exact model and to its approximation.

\subsection{Conditioned Gaussian process}
\label{ssec:cGp}
In this framework, the function $f$ is interpreted as a realization of a Gaussian Process. 
Let us denote $F_\lambda$ a Gaussian process defined, for any $x=(t,\psi)$, as
\begin{equation}
 \label{krigeage}
F_{\lambda}(x)=\sum_{j =1}^L\beta_jh_j(x)+\zeta(x)= \,^tH(x)\bbeta+\zeta(x)\,,
\end{equation}
where $h_1,\ldots,h_L$ are regression functions, $H=(h_1, \ldots, h_L)$, $\bbeta=(\beta_1,\ldots,\beta_L)$ is a vector of parameters, $\zeta$ is a centered Gaussian process 
with covariance function 
$$\text{Cov}(\zeta(x),\zeta(x'))=\sigma^2 K^{\phi}(x,x')\,,$$
where $K^{\phi}$ is a correlation function depending on some parameters $\phi$ and $\lambda=(\bbeta,\sigma,\phi)$ is the vector of all unknown parameters.
For instance, the so-called Gaussian kernel is defined by 
$K^{\phi}(x,x')=\exp(-\phi \Vert x-x' \Vert^2)\,.$ 
The regression functions $h_1,\ldots,h_L$ 
%(and their number $L$)
are usually  
linear functions or low degree polynomials.
The   kernel $K^{\phi}$ has to be chosen with respect to the assumed regularity of the function $f$. Similarly, the regression functions $H$ have to be chosen with respect  to the supposed trend in the function $f$
if some insights on the function are available.
See \cite{koehler:owen:1996,fang:2005} for detailed discussions on the choice of the
regression functions and kernels. 

\medskip

We assume that we are able, in a pre-computation step, to evaluate precisely the function $f$ $n_D$ times for a  
given design of numerical experiments, $D=\{x_1,\ldots x_{n_D}\}$. These "exact" evaluations are denoted  
$z_1=f(x_1),\ldots,z_{n}=f(x_{n_D})$. The $(z_k)$ are different from the $(y_{ij})$ considered in model 
(\ref{eq:MixedModel}) which are noisy observations of $f$ in some unknown points $x=(t, \bpsi)$.
The design of experiments $D$ is usually chosen with respect to a space-filling criterion \citep{fang:2005} in a bounded domain
where the points $(t_{ij},  \psi_i)_{i,j}$ are assumed to be.
For a given kernel $K$, and a given vector $H$ of regression functions, the vector of parameters $\lambda$ has to be estimated. 
Usually, $\lambda$ is  estimated by maximizing the log-likelihood $\ell_F$ of the Gaussian process $F$ (\ref{krigeage}):
\begin{eqnarray}
 \label{eq:vraiskrig}
\ell_F(\lambda; \bz_D)&=&-\frac{1}{2}\log((2\pi\sigma^2)^{n_D}|\Sigma_{{D}{D}}^\phi|)\\
\nonumber
&&-\frac{1}{2}  \,^t(\bz_{D}-H_{D}\bbeta)(\sigma^2\Sigma_{{D}{D}}^\phi)^{-1}(\bz_{D}-H_{D}\bbeta)\,,
 \end{eqnarray}
where 
$(\Sigma_{{D}{D}}^\phi)_{1\leq k,j \leq n_D} =(K^{\phi}(x_k,x_j))$, $(H_{D})_{1\le k\le n_D,1\le j\le L}=h_j(x_k)$. 
The Matlab toolbox DACE \citep{lophaven:etal:2002} provides an optimization algorithm to estimate directly $\lambda$.
We denote by $\hat \lambda =(\hat \bbeta,\hat {\sigma}, \hat \phi)$ the estimates.
The Gaussian process is chosen to be
$F=F_{\hat \lambda}$. It corresponds to a plug-in approach since from now, these parameters are considered as known.

\medskip
Then $f$ is not directly approximated by $F$, but rather by the conditional process denoted $F^D$, defined as the process $F$ conditionally to $F(x_1)=z_1,\ldots,F(x_{n_D})=z_{n_D}$,
in short $\bZ_{D}=\bz_{D}$.
The process $F^{D}$ is still a Gaussian process, defined by its mean and covariance functions, which can be exactly computed. Let us introduce  the partial functions, for any $x\in\mathbb{R}^{d+1}$, 
$K^{\phi}_x:\mathbb{R}^{d+1}\rightarrow\mathbb{R}$  defined by $K^{\phi}_x(x')=K^{\phi}(x,x')$ for any $x'$ and the vector $\Sigma_{x {D}}^{\hat \phi}=\left(K_x^{\hat \phi}(x_k)\right)_{1\le k \le n_D}
$. 
Then, the mean $m_D(x)$ and covariance $C_D(x,x')$ of the process $F^D$  are defined, for all $x, x'$ as
\begin{equation}\label{eq:mean}
m_{D}(x)=H(x)^t\hat\bbeta + \,^t\Sigma_{x {D}}^{\hat \phi}(\Sigma_{{D}{D}}^{\hat \phi})^{-1}(\bz_{D}-H_{D}\hat\bbeta)\,, 
\end{equation}
\begin{equation}
\label{eq:covariance}
C_{D}(x,x') = \hat\sigma^2(K^{\hat \phi}_x(x')
- \,^t\Sigma_{x {D}}^{\hat \phi}(\Sigma_{{D}{D}}^{\hat \phi})^{-1}\Sigma_{x' {D}}^{\hat \phi})\,.
\end{equation}
The mean function $m_{D}$ provides an approximation of the function $f$ for any $x$ and
the variance function $x\mapsto C_{D}(x,x)$ measures the confidence in the accuracy of this approximation. 
The plug-in approach for $\lambda=(\beta,\sigma^2,\phi)$ may lead to underestimate the uncertainty on the quality of the approximation
which is showed to be asymptotically negligible \citep{prasad:rao:1990}.
It can be
used only for the parameter of the correlation kernel $\phi$. In this case, 
when the process is not conditioned on $(\hat \bbeta,\hat {\sigma})$, the conditioned process is a
Student T-process with still closed-form location and scale \citep{santner:etal:2003}. However, we prefer
the complete plug-in approach to   deal with a Gaussian distribution which is easier to incorporate in the mixed model.

\subsection{Kernel interpolation}
We can   interpret the   previous meta-model approximation as a kernel interpolation. Indeed, a Reproducing Kernel Hilbert Space (RKHS) can be associated to $K^{\phi}$ as soon as the kernel $K^{\phi}$ is positive definite. 
% $$\forall (\lambda_1,x_1)\ldots (\lambda_M,x_M)\in\mathbb{R}
%\times \mathbb{R}^{d+1}, \sum_{1\le l,m\le M}\lambda_l\lambda_mK(x_l,x_m)\ge 0\,.$$
Then the partial functions $K^{\phi}_x$ defined in Section \ref{ssec:cGp} span a pre-Hilbert space with inner product  $<K^{\phi}_x,K^{\phi}_{x'}>=K^{\phi}(x,x')$.
Aronszajn's theorem states that there exists a unique completion $\mathcal{H}_K$ of this space with the reproducing property:
\begin{equation*}
\label{reprorkhsmd}
\forall v\in\mathcal{H}_{K},\ x\in\mathbb{R}^{d+1},\ v(x)=<v,K^{\phi}_x> \,.
\end{equation*}
The space $\mathcal{H}_K$ is the  RKHS associated to kernel $K$.  
Then, we focus on the function $g(x)=f(x)-\,^tH(x)\hat \bbeta$, where $\hat \bbeta$ is estimated as before. 
Note that this function can be seen as the residuals of the linear model (\ref{krigeage}) of the random variables $Z_1, \ldots, Z_{n_D}$
on the design $D=\{x_1,\ldots x_{n_D}\}$:
$$ Z_i = f(x_i)= \,^tH(x_i)\hat\bbeta+g(x_i).$$

We consider the following problem of seeking for the function in $\mathcal{H}_K$ which interpolates $g$ on points of $D$ with a minimal norm:
\begin{center}
$\left\{
\begin{array}{l}
\min_{v\in {\cal H}_K}\Vert v \Vert_{{\cal H}_K}\\
g(x_k)=v(x_k), \ k=1,\dots n_D.
\end{array}
\right.$
\end{center}
The solution to this problem is the orthogonal projection of $g$ on the subspace spanned by $(K_{x_1},\ldots,K_{x_{n_D}})$,  denoted  $s_{K,D}(g)$.
If we assume that the function $g=f-\,^tH\hat \bbeta$ belongs to  $\mathcal{H}_K$, then $s_{K,D}(g)$ is defined as 
\begin{eqnarray}
\label{eq:Kinterp}
s_{K,D}(g(x))&=&\,^t\Sigma_{x {D}}^{\hat\phi}(\Sigma_{{D}{D}}^{\hat\phi})^{-1}(\bz_D-H_D\hat \bbeta)\\
&=&\,^t\bu(x)(\bz_D-H_D\hat \bbeta)=\sum_{k=1}^{n_D}u_k(x)g(x_k) \,,\nonumber
\end{eqnarray}
with $\,^t\bu(x)=\,^t\Sigma_{x {D}}^{\hat\phi}(\Sigma_{{D}{D}}^{\hat\phi})^{-1}$ (same notations as in Section \ref{ssec:cGp}). 
This provides an approximation of $f$ which is the same than the function $m_D$ (\ref{eq:mean}).
 The approximation $m_D$ belongs to the RKHS assuming that the regression functions $H$ belong to the RKHS, which is   true for linear or low degree polynomial regressors $H$.

\medskip
The kernel interpolation framework yields to   an upper bound on the point-wise error of this approximation  using 
the reproducing property and a Cauchy-Schwarz inequality.
For any $x$, we have
\begin{eqnarray}
 \label{eq:bornedeterministe}
 \nonumber
 |f(x)-m_D(x)|%&=& |f(x)-\,^tH(x)\hat \bbeta -s_{K,D}(g(x))|\\ 
 &=&|g(x)-s_{K,D}(g(x))|\\
 &\le& \; |<g,K^{\hat\phi}_x-\sum_{k=1}^{n_D}u_k(x)K^{\hat\phi}_{x_k}>|\nonumber\\
 &\le&\Vert g\Vert_{\mathcal{H}_K}\cdot \Vert K^{\hat\phi}_x-\sum_{k=1}^{n_D}u_k(x)K^{\hat\phi}_{x_k}\Vert_{\mathcal{H}_K}\nonumber\\
 &=:& \Vert g\Vert_{\mathcal{H}_K} P_D(x)\,.
\end{eqnarray}
The norm $\Vert g\Vert_{\mathcal{H}_K}$ is unknown and depends on $f$.
The norm $P_D(x)$ does not depend on $f$ (or $g$) but on
 the design of experiments $D$ only.
It holds that
\begin{eqnarray*}\label{eq:borne}
P_D(x)&=&(K^{\hat \phi}(x,x)
-\, ^t\Sigma_{x {D}}^{\hat\phi}(\Sigma_{{D}{D}}^{\hat\phi})^{-1}\Sigma^{\hat\phi}_{x {D}})=\frac 1{\hat\sigma^2} C_D(x,x)\,.
\end{eqnarray*}
Again, there exists a link with the Gaussian process framework: up to the parameter $\hat \sigma^2$,
we obtain the variance function of the conditioned Gaussian process (\ref{eq:covariance}).
\smallskip

For some usual kernels, a uniform  upper-bound on $P_D(x)$ is available as a function of 
$$a_D=\sup_{x\in \mathcal{X}^{d+1}}\min_{1\le k \le n_D}\Vert x-x_k\Vert$$
where $\mathcal{X}$ is a bounded subspace of $\mathcal{R}$.
The value of $a_D$ is related to the coverage of the space $\mathcal{X}$ by the design of experiments. A design of experiments which minimizes this quantity
is said to be minimax \citep{johnson:etal:1990}. The point-wise upper bound is given in the following Proposition \citep{schaback:1995}:
\begin{proposition}\label{prop:errorMetaModel}
Assume that the experimental design $D$ is minimax in $\mathcal{X}$.  Let $\mathcal{H}_K$ denote the RKHS associated to the kernel $K$ which
is assumed to be derived from a radial basis function as proposed in \citet{wu:1992}.
Assume that $f$ lies in $\mathcal{H}_K$.
Let $m_D$ denote the kernel approximation of the function $f$ obtained on the design $D$. 
Then the point-wise error $|f(x)-m_D(x)|$ is uniformly upper-bounded in $\mathcal{X}$ by
$$|f(x)-m_D(x)|\leq \Vert g\Vert_{\mathcal{H}_K} P_D(x)\le \Vert g\Vert_{\mathcal{H}_K} G_K(a_D)\,.$$
where the function $G_K$ is defined on $\mathbb{R}^+$ and is such that $\underset{a\rightarrow 0^+}\lim G_K(a)=0$.

Furthermore, if the regressors $H\in \mathcal{H}_K$, then there exists a constant $C>0$ such that
$$|f(x)-m_D(x)|\leq  C \Vert f\Vert_{\mathcal{H}_K} G_K(a_D)\,.$$
\end{proposition}
%If we assume that the regression functions are in the RKHS ${\mathcal{H}_K}$, this inequality still holds up to a constant if we replace $\Vert g\Vert_{\mathcal{H}_K}$ by $\Vert f\Vert_{\mathcal{H}_K}$.
For instance, when using a Gaussian kernel $K^\phi(x,x')=e^{-\phi\Vert x-x'\Vert^2}$, the function $G_K$ is
$G_K(a)=Ce^{-\delta/ a^2}$ where $C$ and $\delta$ are constants depending on $\phi$.

%%%%%%%%%%%%% 
\section{Three mixed meta-models}\label{sec:MixedMetaModel}
The estimation of the population parameter $\theta$  is performed on the meta-model approximation of the mixed model (\ref{eq:MixedModel}). For computational reasons, we introduce three mixed meta-models. 

\subsection{Complete mixed meta-model}
Let us introduce  the so-called complete mixed meta-model that integrates the meta-model approximation. 
The regression function $f$ in (\ref{eq:MixedModel}) is approximated by $  F_D(t, \psi) = m_D(t,\psi)+r(t, \psi)$:
\begin{eqnarray} \label{eq:completeMixedMetaModel}
y_{ij} &=&  F_{D}(t_{ij}, \psi_i)  + \sigma_\varepsilon\, \varepsilon_{ij},  \ \varepsilon_{ij} \sim_{iid} \mathcal{N}(0,1), \\
 \psi_i &\sim_{iid}&    \mathcal{N}(\mu,\Omega), \nonumber\\
 F_D(t,\psi) &=&   m_{D}(t,\psi)+ r(t,\psi) , \quad \mbox{with,} \nonumber\\
 r(t, \psi) &\sim& \mathcal{GP} ( 0, C_{D}(t,\psi; t,\psi) ).\nonumber
\end{eqnarray}
  Whereas model (\ref{eq:MixedModel}) is homoscedastic (constant error variance), the mixed metamodel (\ref{eq:completeMixedMetaModel}) is heteroscedastic. 
Let us emphasize   that \textit{this is not a standard heteroscedastic error model}. Indeed, we have:
 $$y_{ij}|\psi_i \sim \mathcal{N}(m_{D}(t_{ij}, \psi_i), \Gamma_D(t_{ij}, \psi_i))$$
  with $$\Gamma_D(t,\psi) = \sigma_\varepsilon^2 + C_{D}(t,\psi; t,\psi)\,,$$ but the $(y_{ij}|\psi_i)_j$ are \textit{not independent}, as well as the individual vectors $(\by_i)_i$. This is due to the fact that the $ (r(t_{ij},\psi_i))_{ij}$ are realizations of the same Gaussian process. 
  This is a major difference with approximations that have already been  proposed in the literature \citep{Donnet2007, Grenier2014}. Especially, this complicates the implementation of the MCMC scheme. %But we expect to gain in the accuracy of the estimator with this approach. 

 \medskip
 We propose to estimate $\theta$ as the maximum of the likelihood of  model (\ref{eq:completeMixedMetaModel}). 
 We denote $$\bm_D(\bt,\bpsi) = (m_D(t_{ij},\psi_i))_{1\leq i\leq N, 1\leq j\leq n_i}$$ the vector of the approximate mean, evaluated on 
 $(\bt, \bpsi) = (t_{ij}, \psi_i)_{1\leq i\leq N, 1\leq j\leq n_i}$. Similarly, we denote 
 $$\bC_D(\bt,\bpsi) = (C_{D}(t_{ij},\psi_i; t_{i'j'},\psi_{i'}))_{1\leq i,i'\leq N, 1\leq j,j'\leq n_i}\,.$$
 The likelihood of model (\ref{eq:completeMixedMetaModel}) is then:
\begin{eqnarray}\label{eq:completeMetaLikelihood}
p_D(\by;\btheta)&=&
%&&=\int p(\bpsi;\btheta) p_D(\by|\bpsi;\btheta)d\bpsi\,,\nonumber\\
\int \Bigg\{ p(\bpsi;\btheta)\;  \frac{1}{(2\pi)^{n_{tot}/2}|\sigma_\varepsilon^2\,I_{n_{tot}}+\bC_D(\bt,\bpsi)|^{1/2}}\\
&&\exp\bigg(-\frac{1}{2}\,^t(\by-\bm_D(\bt,\bpsi))(\sigma_\varepsilon^2\, I_{n_{tot}}+\bC_D(\bt,\bpsi))^{-1}(\by-\bm_D(\bt,\bpsi)) \bigg)\,d\bpsi
\Bigg\}\nonumber.
\end{eqnarray}
This likelihood is not explicit because   function $m_{D}(t_{ij}, \psi_i)$ is not linear in $\psi_i$.  
As said previously, this likelihood  cannot be simplified as a product of individual likelihoods because the $\by_i$ are not independent 
(the matrix $\bC_D(\bt,\bpsi)$ is a full matrix). 
The corresponding estimation algorithm   requires to invert this $n_{tot}\times n_{tot}$-matrix at each iteration (at least $N\times 2d$ per iteration), which  is highly computationally intensive. Therefore, we  introduce an intermediate mixed meta-model 
by considering only the diagonal of $\bC_D(\bt,\bpsi)$. 

%\textit{a enlever ?}
%Furthermore, the likelihood of the complete data (the integrand) does not belong to the 
%exponential family with respect to the parameter $\sigma_\varepsilon$ because of the heteroscedasticity of the model. \textbf{ c'est faux, il suffit de considerer $r$ dans les donnees completes}... 
%As this is an assumption for running the SAEM algorithm, $\sigma_\varepsilon$ might be more difficult to estimate. 

\subsection{Intermediate mixed meta-model}

In the intermediate mixed meta-model, the regression function $f$ is approximated by $ m_D(t,\psi)+\bar{r}(t, \psi)$, where $\bar{r}(t, \psi)$ has a diagonal covariance matrix $\Lambda_{i, \psi_i}=\diag (C_{D}(\bt_i,\psi_i) )$:
\begin{eqnarray} \label{eq:intermediateMixedMetaModel}
y_{ij} &=&  m_{D}(t_{ij}, \psi_i) + \bar{r}(t_{ij}, \psi_i)  + \sigma_\varepsilon\, \varepsilon_{ij}, \\
\varepsilon_{ij} &\sim_{iid}& \mathcal{N}(0,1) \,, \nonumber\\
 \psi_i &\sim_{iid}&    \mathcal{N}(\mu,\Omega) \nonumber\\
\bar{r}(\bt_i,\psi_i) &\sim_{ind}& \mathcal{GP} ( 0, \Lambda_{i, \psi_i}=\diag (C_{D}(\bt_i,\psi_i) )).\nonumber
\end{eqnarray}
The likelihood of  model  (\ref{eq:intermediateMixedMetaModel}) is then:
\begin{eqnarray}\label{eq:intermediateMetaLikelihood}
\bar{p}_D(\by;\btheta)
&=&\prod_{i=1}^N \int\Bigg\{ p(\psi_i;\btheta)
\frac{1}{(2\pi)^{n_i/2}  \prod_{j=1}^{n_i}(\sigma_\varepsilon^2+C_D((t_{ij},\psi_i),(t_{ij},\psi_i)))^{1/2}} \\
&&\exp\bigg(-\frac{1}{2}\,^t(\by_i-m_D(\bt_i,\psi_i))(\sigma_\varepsilon^2 I_{n_i}+ \Lambda_{i, \psi_i})^{-1}(\by_i-m_D(\bt_i,\psi_i)) \bigg)
 d\psi_i \Bigg\}\,.\nonumber
\end{eqnarray}
This form of the likelihood is separable with respect to $\psi_i$ and can be written as a product over the individuals which are independent.  The covariance matrix $\sigma_\varepsilon^2 I_{n_i}+ \Lambda_{i, \psi_i}$ is diagonal and can be easily inverted. This will substantially reduce the computation time of the estimation method. However the intermediate model is heteroscedastic, and $\sigma_\varepsilon$ might   be more difficult to estimate than in the exact model. This is why we introduce a simpler mixed meta-model.

\subsection{Simple mixed meta-model}
The simple mixed meta-model neglects the error of approximation of the computer model. The regression function is then $m_D$:
\begin{eqnarray} \label{eq:simpleMixedMetaModel}
y_{ij} &=&  m_{D}(t_{ij}, \psi_i) + \sigma_\varepsilon\, \varepsilon_{ij},  \quad \varepsilon_{ij} \sim_{iid} \mathcal{N}(0,1) ,\\
 \psi_i &\sim_{iid}&    \mathcal{N}(\mu,\Omega) .\nonumber
\end{eqnarray}

The simple mixed meta-model (\ref{eq:simpleMixedMetaModel}) has similar properties than  model (\ref{eq:MixedModel}): it is  homoscedastic (constant error variance), \textit{the vectors $(\by_i)_i$ are independent and identically distributed}, and for each individual $i$, conditionally to $\psi_i$, the $(y_{ij})_j$ are independent. 
The likelihood of model (\ref{eq:simpleMixedMetaModel}) is given by: 
\begin{eqnarray}\label{eq:simpleMetaLikelihood}
\tilde p_D(\by;\btheta)&=&\prod_{i=1}^N \int \Bigg\{p(\psi_i;\btheta)
\frac{1}{(2\pi \sigma_\varepsilon^2)^{n_i/2}} \\
&&
\exp\bigg(-\frac{1}{2}\,^t(\by_i-m_D(\bt_i,\psi_i))(\sigma_\varepsilon^2 I_{n_i})^{-1}
(\by_i-m_D(\bt_i,\psi_i)) \bigg)
 d\psi_i \Bigg\}, \nonumber
\end{eqnarray}
which has the same form than likelihood of model (\ref{eq:MixedModel}).

%%%%%%%%%%%%%  
\section{Population parameter estimation}\label{sec:estimation}
Likelihoods of the mixed meta-models being not explicit, we    resort to the family of EM algorithm to estimate the parameters $\theta$, which is a classical approach for models with non-observed or incomplete data. 
 We start with the SAEM algorithm for the exact mixed model and then for the three mixed meta-models.

\subsection{Estimation for the exact mixed model}
The objective is to maximize the likelihood $p(\by; \theta)$ of the exact mixed model (\ref{eq:MixedModel}). 
Let us briefly cover the EM principle  \citep{Dempster1977}. 
 The complete data of the mixed  model is  $(\by,\bpsi)$.
The EM algorithm maximizes the
 $Q(\theta|\theta')=\mathbb{E}(L(\by,\bpsi;\theta)|\by;\theta')$ function
in 2 steps, where $L(\by,\bpsi;\theta)$ is the log-likelihood of the complete data for the   mixed  model (\ref{eq:MixedModel}) and $\mathbb{E}$ is the   expectation under the conditional distribution $p(\bpsi|\by; \theta')$. % $(y,\psi)$.
At the $k$-th iteration, the E step is the
evaluation of  $Q_k(\theta)=Q(\theta \, \vert \,\widehat{\theta}^{(k-1)})$, whereas the M step
updates $\widehat{\theta}^{(k-1)}$ by maximizing $Q_{k}(\theta)$.
For cases with a non analytic E step,   
\cite{Delyon1999}
introduce a stochastic version SAEM of the EM algorithm which
evaluates the integral $Q_k(\theta)$ by a stochastic approximation procedure.
The E step is then divided into a
simulation step (S step) of the missing data $\bpsi^{(k)}$ under the
conditional distribution $p(\bpsi|\by;  \widehat{\theta}^{(k-1)})$ and a stochastic
approximation step (SA step) of the conditional expectation, using $(\gamma_k)_{k\geq 0}$ a sequence of positive numbers decreasing to 0: 
$$
Q_{k}(\theta)=Q_{k-1}(\theta)+\gamma_k( L(\by,\bpsi^{(k)}; \theta)-Q_{k-1}(\theta)).
$$
In cases where the simulation of the non-observed vector $\psi$ cannot be directly performed,   
\cite{Kuhn2005} propose to combine the SAEM algorithm with a Markov Chain Monte Carlo (MCMC) procedure. 
The idea is to simulate a Markov chain   $\bpsi^{(k)}$ 	by use of  a Metropolis-Hastings (M-H) algorithm 
with $p(\bpsi|\by; \widehat{\theta}^{(k-1)})$  as the unique stationary distribution.

The complete data likelihood $L(\by,\bpsi;\theta)$ of the exact mixed model belongs to the regular curved exponential 
family:
 $$p(\by, \bpsi;\theta)= \exp \left\{-\nu(\theta)+ <S(\bpsi), \lambda(\theta)>\right\}$$
where $<\cdot, \cdot>$ denotes the scalar product, the minimal sufficient statistic $S(\by,\bpsi)$  takes its values in an open subset $\mathcal{S}$ of $\mathbb{R}^m$, $\nu_D$ and $\lambda$ are functions of $\theta$.  
Then
%\begin{itemize}
%\item[\textbf{(M1)}] $L_D(y,\psi;\theta)= -\lambda_D(\theta)+\langle S_D(y,\psi),\nu_D(\theta)\rangle$ ,
%\end{itemize}
%where $\left\langle .,.\right\rangle$ is the scalar product, $\lambda_D$ and $\nu_D$ are  functions of $\theta$ and  $S_D(y,\psi)$ is  the minimal sufficient statistic of the complete data likelihood.  
the SA step reduces to approximate  $\mathbb{E}\left\lbrack  S(\by,\bpsi) \vert  \widehat{\theta}^{(k-1)} \right\rbrack$ at each iteration by the value $s_k$. 
The sufficient statistics for the exact mixed model are   classically $S_{1}(\by, \bpsi) = \sum_{i=1}^N \psi_i$, $S_{2}(\by, \bpsi) = \sum_{i=1}^N \psi_i \, ^t\psi_i$ and\\
$S_{3}(\by, \bpsi) = \sum_{i=1}^N\sum_{j=1}^{n_i} (y_{ij}-f(t_{ij}, \psi_i))^2$  \citep{Samson2007}. 
Then the M step is explicit and easy to implement. The convergence of the SAEM-MCMC algorithm has been proved when the complete data likelihood  belongs to the regular curved exponential family and under additional assumptions (see Proposition \ref{prop:Kuhn}). Thus the exponential family plays a crucial role to obtain an efficient algorithm. 

%We will see that this is not always the case, especially for the complete and the intermediate mixed meta-models. 
 
\subsection{Estimation for the simple mixed meta-model}
%As explained before, also the complete and intermediate mixed meta models have the advantage to directly incorporate the approximation error, their implementation is tricky, especially to estimate $\sigma_\varepsilon$. We now present how the implementation of the SAEM-MCMC is clearly easier when working with the simple mixed meta-model. 
 
The objective is to maximize the likelihood $\tilde p_D(\by; \theta)$ of the simple mixed meta-model (\ref{eq:simpleMixedMetaModel}). In the following, all the quantities referring to this approximate likelihood $\tilde p_D(\by; \theta)$ are indexed by $D$ with a tilde symbol. 
The corresponding complete data likelihood 
$\tilde L_D(\by,\bpsi;\theta)$ belongs to the regular curved exponential family with minimal sufficient statistics  $\tilde S_D(\by,\bpsi)$, which are the same as the exact mixed model.  
In that model, the MCMC algorithm is easy to implement because of the independence of the observations of the individuals. %Furthermore, the maximization step with respect to $\sigma_\varepsilon$ is explicit. 
More precisely, the SAEM-MCMC is described as follows.

\begin{algo}\textbf{(SAEM-MCMC algorithm for the simple mixed meta-model)}\\\label{algo2}
At iteration $k$, given the current values of the estimators  $\hat\mu^{(k-1)}, \hat\Omega^{(k-1)},\hat \sigma_\varepsilon^{2\,(k-1)}$: %, the SAEM algorithm proceeds as follows:
\begin{itemize}
\item[] \textbf{Simulation step:} For each individual $i$ separately and successively, update   $\psi_i^{(k)}$  with $m$ iterations of an MCMC procedure with $\tilde p_D(\psi_i|\by_i; \widehat{\theta}^{(k-1)})$ as stationary distribution:

 \noindent For $l=1\ldots,m$, given a current value $\psi_i^{l-1}$   for individual $i$:
	\begin{itemize}
	\item Simulate a candidate $\psi_i^c$  with a proposal distribution $q_{\widehat{\theta}^{(k-1)}}(\cdot| \psi_i^{l-1})$. 
	\item \textbf{Meta-model step:}   Evaluate the meta-model\\ $m_{D}(t_{ij}, \psi_i^{c})$ for all $j=1, \ldots, n_i$. 
%	\item Compute the acceptation probability $\tilde\alpha_i(\psi_i^c, \psi_i^{l-1})$
%$$\tilde\alpha_i(\psi_i^c, \psi_i^{l-1})= \min \left(\frac{\tilde p_D(\by_i|\psi_i^c; \widehat{\theta}^{(k-1)}) p(\psi_i^c;\widehat{\theta}^{(k-1)}) }{\tilde p_D(\by_i|\psi_i^{l-1}; \widehat{\theta}^{(k-1)})p(\psi_i^{l-1};\widehat{\theta}^{(k-1)})}\frac{q_{\widehat{\theta}^{(k-1)}}(\psi_i^c|\psi_i^{l-1})}{q_{\widehat{\theta}^{k-1}}(\psi_i^{l-1}|\psi_i^c)}, 1\right)
%$$

	\item The candidate is accepted,  $\psi_i^{l}=\psi_i^c$, with   probability $\tilde\alpha_i(\psi_i^c, \psi_i^{l-1})$; otherwise the candidate is rejected, $\psi_i^{l}=\psi_i^{l-1}$  with   probability $1- \tilde\alpha_i(\psi_i^c, \psi_i^{l-1})$, where
	\begin{eqnarray*}
	\tilde\alpha_i(\psi_i^c, \psi_i^{l-1})
	&=& \min \bigg(\frac{\tilde p_D(\by_i|\psi_i^c; \widehat{\theta}^{(k-1)}) p(\psi_i^c;\widehat{\theta}^{(k-1)}) }{\tilde p_D(\by_i|\psi_i^{l-1}; \widehat{\theta}^{(k-1)})p(\psi_i^{l-1};
	\widehat{\theta}^{(k-1)})}
	\frac{q_{\widehat{\theta}^{(k-1)}}(\psi_i^{l-1}|\psi_i^c)}{q_{\widehat{\theta}^{(k-1)}}(\psi_i^c|\psi_i^{l-1})}, 1\bigg)\nonumber\,. 
	\end{eqnarray*}

	\end{itemize}
\noindent Set $\psi_i^{(k)}=\psi_i^{m}$.
\item[] \textbf{Stochastic Approximation step:} update   the sufficient statistics: %$\tilde S_D$
	\begin{eqnarray*}
	s_{k,1} & = & s_{k-1,1} + \gamma_k\,\left(\sum_{i=1}^N \psi_i^{(k)} -s_{k-1,1}\right)\nonumber \,,\\
	s_{k,2} & = & s_{k-1,2} + \gamma_k\,\left(\sum_{i=1}^N \psi_i^{(k)}\,^t\psi_i\,^{(k)}\, -s_{k-1,2}\right)\nonumber \,,\\
	s_{k,3} & = & s_{k-1,3} 
	+ \gamma_k\,\left(\sum_{i=1}^N\sum_{j=1}^{n_i}(y_{ij} - m_{D}(t_{ij}, \psi_i^{(k)}))^2-s_{k-1,3}\right)\nonumber\,. 
	\end{eqnarray*}
\item[] \textbf{Maximisation step:} update the population parameters
	\begin{eqnarray*}
	\widehat{ \mu}^{(k)} & = & \frac{s_{k,1}}{N}, \quad	 \widehat{\Omega}^{(k)}  =  \frac{s_{k,2}}{N}-\frac{s_{k,1}\, ^ts_{k,1}}{N^2}, \quad 
	\widehat{\sigma_\varepsilon}^{2\,^{(k)}}  =  \frac{s_{k,3}}{n_{tot}}\,.
	\end{eqnarray*}
		\end{itemize} 
\end{algo}

\subsection{Estimation for the intermediate mixed meta-model} 
In the following, all the quantities referring to the approximate likelihood $\bar p_D(\by; \theta)$ of the intermediate mixed meta-model are indexed by $D$ with a bar symbol. 
%This intermediate mixed meta-model has two advantages with respect to the   other mixed models: it is conservative since it   takes into account the meta-model approximation of the meta-model, but it remains reasonable because the vectors $\by_i$ are independent. 

This model belongs to the exponential family when  the Gaussian process $\bar{r}$ is considered in the hidden states. Then the complete data of the intermediate mixed meta-model are    $(\by,\bpsi, \bar \br)$ where $\bar\br=(\bar r(t_{ij},\psi_i))_{i=1, \ldots, N, j=1, \ldots, n_i}$.
The complete log-likelihood is thus: 
\begin{eqnarray*}
 \bar L_D(\by,\bpsi, \bar\br;\theta)
 &=& \log \bar p_D(\by|\bar\br, \bpsi; \theta) + \log \bar p_D(\bar\br|\bpsi; \theta) + \log p(\bpsi;\theta)\\
 &=& cst  -  \frac{n_{tot}}{2} \log(\sigma_\varepsilon^2)
 -\frac12 \sum_{ij} \frac{(y_{ij}-m_D(t_{ij}, \psi_i)-\bar r(t_{ij},\psi_i))^2}{\sigma_\varepsilon^2}\\
&&-\frac12\sum_i\log(| \Lambda_{i, \psi_i}|)
 -\frac12\sum_i  \,^t\bar\br \Lambda_{i, \psi_i}^{-1}\bar \br -\frac N2 \log(|\Omega|)\\
&& -\frac12\sum_i \,^t(\psi_i-\mu)\Omega^{-1}(\psi_i-\mu)\,,
\end{eqnarray*}
where $cst$ denotes a constant term.
The E-step is the computation of  
\begin{eqnarray*}
Q(\theta|\widehat{\theta}^{(k-1)})
&=&\mathbb{E}(\bar L_D(\by,\bpsi, \bar \br;\theta)|\by;\widehat{\theta}^{(k-1)})\\
&=&
 \int \int \log \bar p_D(\by, \bar\br, \bpsi; \theta) \bar p_D(\bar\br, \bpsi|\by; \widehat{\theta}^{(k-1)}) d\bar\br d\bpsi\\
&=& 
\int\left( \int \log \bar p_D(\by,\bar\br, \bpsi; \theta) \bar p_D(\bar\br|\by, \bpsi; \widehat{\theta}^{(k-1)}) d\bar\br\right)
\bar p_D(\bpsi|\by; \widehat{\theta}^{(k-1)})  d\bpsi\,.
\end{eqnarray*}
The conditional distribution $ \bar p_D(\bar\br_i|\by_i, \psi_i; \widehat{\theta}^{(k-1)})$ is explicit, Gaussian, with mean and covariance defined by 
\begin{eqnarray*}
\bar m_{r,\psi_i}^{(k-1)}  &=&\bar \Gamma_{r,\psi_i}^{(k-1)}  (\by_i-m_D(\psi_i))/\widehat{\sigma_\varepsilon^{2}}^{(k-1)} \,,\\
 \bar \Gamma_{r,\psi_i}^{(k-1)} & =&  (1/\widehat{\sigma_\varepsilon^{2}}^{(k-1)}+ \Lambda_{i, \psi_i}^{-1})^{-1}.  
\end{eqnarray*}
Integrated with respect to $\bar{\br}$ inside $Q(\theta|\widehat{\theta}^{(k-1)})$ yields to
\begin{eqnarray*}
Q(\theta|\widehat{\theta}^{(k-1)})
&=& \int  \bigg[- \frac{n_{tot}}2 \log(\sigma_\varepsilon^2)  
-\frac12\sum_{i} \frac{\|\by_i-m_D(\bt_i, \psi_i)-\bar m_{r, \psi_i}^{(k-1)}\|^2}{\sigma_\varepsilon^2} \\
&&-\frac1{2\sigma_\varepsilon^2} \sum_i \Tr(\bar \Gamma_{r, \psi_i}^{(k-1)})-\frac12\sum_i \log | \Lambda_{i, \psi_i}|\\
&&-\frac12 \sum_i  \,^t\bar  m_{r, \psi_i}^{(k-1)} \Lambda_{i, \psi_i}^{-1}\bar m_{r, \psi_i}^{(k-1)}-\frac12 \sum_i \Tr( \Lambda_{i, \psi_i}^{-1}\bar \Gamma_{r, \psi_i}^{(k-1)})\\
&&
 -\frac N2\log(|\Omega|) -\frac12\sum_i \,^t(\psi_i-\mu)\Omega^{-1}(\psi_i-\mu)
\bigg]\\
&&p(\bpsi|\by; \widehat{\theta}^{(k-1)})   d\bpsi+cst\,.
\end{eqnarray*}

Then the  sufficient statistic corresponding to $ \sigma_\varepsilon^2$ is changed to $\bar S^{(k-1)}_{D,3}(\by, \bpsi, \br) = \sum_{i=1}^N \|\by_{i}-m_D( \psi_i)-\bar m_{r,\psi_i}^{(k-1)}\|^2+  \Tr(\bar \Gamma_{r,\psi_i}^{(k-1)})$. 
The simulation step is a standard one, which can be applied to each individual separately. The MCMC algorithm targets $\bar p_D(\psi_i|\by_i; \widehat{\theta}^{(k-1)})$ as stationary distribution, where the process $\bar \br_i$ has been integrated out. The acceptance probability only   requires the knowledge of $\bar p_D(\by_i|\psi_i; \theta^{(k-1)})$ which is a Gaussian density with covariance matrix $\sigma_\varepsilon^2I_{n_i} +\Lambda_{i, \psi_i}$. As this matrix is diagonal, its inversion   at each iteration is fast. Finally, the SAEM-MCMC proceeds as follows:

\begin{algo}\textbf{(SAEM-MCMC algorithm for the intermediate mixed meta-model)}\\\label{algo1}
At iteration $k$, given the current values of the estimators  $\hat\mu^{(k-1)}, \hat\Omega^{(k-1)},\hat \sigma_\varepsilon^{2\,(k-1)}$: %, the SAEM algorithm proceeds as follows:

\begin{itemize}
\item[] \textbf{S step:} For each individual $i$ separately and successively, update   $\psi_i^{(k)}$  with $m$ iterations of an MCMC procedure with $\bar p_D(\psi_i|\by_i; \widehat{\theta}^{(k-1)})$ as stationary distribution.
 
\item[] \textbf{SA step:} update   the sufficient statistics $s_{k,1}$ and $s_{k,2}$ as usual and update 
	\begin{eqnarray*}
	s_{k,3} & = & s_{k-1,3} + \gamma_k\,\bigg(\sum_{i=1}^N \|\by_{i}-m_D( \psi_i)-\bar m_{r,\psi_i^{(k)}}^{(k-1)} \|^2
	+  \Tr\, \left(\bar \Gamma_{r,\psi_i^{(k)}}^{(k-1)}\right)  -s_{k-1,3}\bigg)\nonumber 
	\end{eqnarray*}
\item[] \textbf{M step:} as usual.

\end{itemize}
\end{algo}

\subsection{Estimation for the complete mixed meta-model}
% We now detail the SAEM-MCMC for the complete mixed meta-model (\ref{eq:completeMixedMetaModel}). 
In the following, all the quantities referring to the approximate likelihood $p_D(\by; \theta)$ of  the complete mixed meta-model (\ref{eq:completeMixedMetaModel}) are indexed by $D$. 

The main difficulty comes from the fact that model (\ref{eq:completeMixedMetaModel}) is heteroscedastic
but not in a standard way: %. Two consequences:  
  the conditional distributions of $\by_i|\psi_i$ are not independent and all the subjects have to be treated together. 
  Similarly as the intermediate model, we consider the Gaussian process $r$ in the complete data and we have to integrate out with respect to $r$ to compute the function $Q$.  The conditional distribution
  $ p_D(\br|\by, \bpsi; \widehat{\theta}^{(k-1)})$ is explicit, Gaussian, with mean and covariance defined by: 
  \begin{eqnarray*}
m_{r, \bpsi}^{(k-1)}  &=& \Gamma_{r, \bpsi}^{(k-1)}  (\by-m_D(\bpsi))/\widehat{\sigma_\varepsilon^{2}}^{(k-1)}, \\
\Gamma_{r, \bpsi}^{(k-1)} & =&  (1/\widehat{\sigma_\varepsilon^{2}}^{(k-1)}+\bC_D(\bt,\bpsi)^{-1})^{-1}.    
  \end{eqnarray*}
The matrix $\Gamma_{r,\bpsi}^{(k-1)}  $ has dimension $n_{tot}\times n_{tot}$ and cannot be split as it was the case with the intermediate model. 
Thus the inversion of $ \Gamma_{r}$ and $C_D$ increases dramatically  the computation time of the estimation algorithm. 

  %, 2/ the complete data likelihood  $L_D(\by,\bpsi;\theta)$ does not belong to the regular curved exponential family. \textbf{la deuxieme set fausse}. 
 
%Let us detail 1/. 
Moreover, the MCMC step is also more complex. 
Indeed, the conditional distributions $p_{D}(\bpsi|\by)$ cannot be written as a product of individual conditional distributions. But the MCMC kernels are applied to each subject $i$ successively. The corresponding target distribution is the conditional distribution $p_{D}(\psi_i|\by, \bpsi_{-i})$
where $\bpsi_{-i}=(\psi_1,\ldots, \psi_{i-1}, \psi_{i+1}, \ldots,\psi_N)$ is the vector of individual parameters except individual $i$ (with obvious notations when $i=1$ 
or $i=N$) and not the distribution $p_{D}(\psi_i|\by_i)$ as in a standard heteroscedastic mixed model. 
This increases the difficulty of implementation of the MCMC: the whole covariance function $\bC_D(\bt, \bpsi)$, evaluated at each point $(t_{ij}, \psi_i)$,
has to be evaluated and inverted at each iteration of the MCMC scheme.  

%Let us detail 2/. Parameters $\mu$ and $\Omega$ can be estimated as for the exact mixed model, with sufficient statistics $S_{1}(\by, \bpsi) = \sum_{i=1}^N \psi_i$, $S_{2}(\by, \bpsi) = \sum_{i=1}^N \psi_i \, ^t\psi_i$. The problem exists for the estimation of the error variance $\sigma_\varepsilon^2$. Indeed, due to the full covariance function $C_D$, there is no sufficient statistic for $\sigma_\varepsilon^2$. We propose to estimate this parameter by a numerical maximization of the complete likelihood during the M-step. This will induce a computational cost. 

\begin{algo}\textbf{(SAEM-MCMC algorithm for the complete mixed meta-model)}\\
\label{algo1}
At iteration $k$, given the current values of the estimators  $\hat\mu^{(k-1)}, \hat\Omega^{(k-1)},\hat \sigma_\varepsilon^{2\,(k-1)}$: %, the SAEM algorithm proceeds as follows:

\begin{itemize}
\item[] \textbf{S step:} for each individual $i$ successively, given the current values 
$\bpsi_{-i}^{(k)}=(\psi_{1}^{(k)}, \ldots,$ $ \psi_{i-1}^{(k)}, \psi_{i+1}^{(k-1)}, \ldots, \psi_{N}^{(k-1)})$ of all the other individuals,
update  $\psi_i^{(k)}$  with $m$ iterations of an MCMC procedure with $p_D(\psi_i|\by, \bpsi_{-i}^{(k)}; \widehat{\theta}^{(k-1)})$ as stationary distribution:

\noindent For $l=1\ldots,m$, given a current value $\psi_i^{l-1}$   for individual $i$ and a current vector\\
$\bpsi^{(k)l-1}=(\psi_{1}^{(k)}, \ldots, \psi_{i-1}^{(k)}, \psi_i^{l-1}, \psi_{i+1}^{(k-1)}, \ldots, \psi_{N}^{(k-1)})$ for all individuals:
	\begin{itemize}
	\item Simulate a candidate $\psi_i^c$  with a proposal distribution $q_{\widehat{\theta}^{(k-1)}}(\cdot| \psi_i^{l-1})$.
	\item Set $\bpsi^c= (\psi_{1}^{(k)}, \ldots, \psi_{i-1}^{(k)},\psi_i^c, \psi_{i+1}^{(k-1)}, \ldots, \psi_{N}^{(k-1)})$. 
	\item \textbf{Meta-model step:}  For all $j=1, \ldots, n_i$, evaluate the meta-model $m_{D}(t_{ij}, \psi_i^{c})$. For all subjects $i',i''=1, \ldots, N$ (including subject $i$) and all observations $j',j''$, evaluate  the covariance functions
	$C_D(t_{i'j'}, \bpsi_{i'}^{c}; t_{i''j''}, \bpsi_{i''}^{c})$  and invert the obtained matrix $C_D$. 
	\item The candidate is accepted,  $\psi_i^{l}=\psi_i^c$, with   probability $\alpha_i(\psi_i^c, \psi_i^{l-1})$; otherwise   $\psi_i^{l}=\psi_i^{l-1}$ with   probability $1- \alpha_i(\psi_i^c, \psi_i^{l-1})$ where 
\begin{eqnarray*}
\alpha(\psi_i^c, \psi_i^{l-1})
&=& \min \Bigg(\frac{p_D(\by|\bpsi^c; \widehat{\theta}^{(k-1)}) p(\psi_i^c;\widehat{\theta}^{(k-1)}) }{p_D(\by|\bpsi^{(k)l-1}; \widehat{\theta}^{(k-1)})p(\psi_i^{(l-1)};\widehat{\theta}^{(k-1)})}
\frac{q_{\widehat{\theta}^{(k-1)}}(\psi_i^{(l-1)}|\psi_i^c)}{q_{\widehat{\theta}^{(k-1)}}(\psi_i^c|\psi_i^{(l-1)})}, 1\Bigg)\,. 
\end{eqnarray*}
	\end{itemize}
\noindent Set $\psi_i^{(k)}=\psi_i^{m}$.
	\item[] \textbf{SA step:} update   the sufficient statistics $s_{k,1}$ and $s_{k,2}$ as before and update: 
	\begin{eqnarray*}
	%s_{k,1} & = & s_{k-1,1} + \gamma_k\,\left(\sum_{i=1}^N \psi_i^{(k)} -s_{k-1,1}\right)\nonumber \\
	%s_{k,2} & = & s_{k-1,2} + \gamma_k\,\left(\sum_{i=1}^N \psi_i^{(k)}\,^t\psi_i\,^{(k)}\, -s_{k-1,2}\right)\nonumber \\
	s_{k,3} & = & s_{k-1,3} + \gamma_k\,\bigg(  \|\by-m_D(\bpsi)-m_{r,\bpsi^{(k)}}^{(k-1)}\|^2
	+  \Tr\, \left(\Gamma_{r,\bpsi^{(k)}}^{(k-1)}\right)  -s_{k-1,3}\bigg)\nonumber\,. 
	\end{eqnarray*}
\item[] \textbf{M step:} as usual.
 
\end{itemize}
\end{algo}

Let us emphasize that the MCMC  in the S step  is difficult to implement due to the heteroscedasticity of the complete mixed meta-model.  
This MCMC algorithm may have poor mixing properties because the vectors $\psi_i$ are updated successively while they are highly correlated through this non-diagonal matrix $\bC_D(\bt, \bpsi)$.
Another solution could be to design a proposal in the MCMC algorithm for the whole vector $\bpsi$. However, such a proposal is quite complicated to construct since
the dimension of $\bpsi$ is high: $d\times N$.

%The   intermediate and simple mixed meta-models will have the advantage to reduce the complexity of the MCMC algorithm (see the two next sections).  

\subsection{Fisher Information matrix estimates}
 
 Using formula in \cite{Louis1982} and estimation scheme proposed in \cite{Delyon1999}, confidence intervals can be obtained on the parameters implementing a stochastic approximation scheme of the Fisher Information matrix. It is only necessary to approximate the gradient and the Hessian matrix of the log-likelihood of the complete data:
 %Some computations on the complete likelihood are necessary.
 \begin{equation}
\log p(\bpsi,\by;\btheta)= \log p(\by|\bpsi;\btheta)+\log p(\bpsi;\btheta) \,.
 \end{equation}
Actually, $\log p(\by|\bpsi;\btheta)$ does not depend on $\mu$ et $\Omega$, hence the gradient and Hessian computations
are only about $\log p(\bpsi; \btheta)$ which is a multivariate normal $\mathcal{N}(\mu,\Omega)$.
Thus, this implementation does not depend on the mixed model and remains the same for the standard mixed model and the three mixed meta-models.

%%%%%%%%%%%%%  
\section{Convergence of the SAEM algorithm to the maximum likelihood of the meta-model}\label{sec:ConvApproxMLE}

Since the SAEM-MCMC algorithm is not applied to  model
(\ref{eq:MixedModel}), but to an approximate mixed model, it is not possible to prove the convergence of the algorithm toward a local maximum of the exact
likelihood $p(\by;\btheta)$. 
However, it is possible to apply the results of \cite{Kuhn2005} 
for the three mixed meta-models. Hence, the algorithms   converge toward a local maximum of the 
likelihood $p_D(\by;\theta)$, $\bar p_D(\by;\theta)$  and $\tilde p_D(\by;\theta)$ when applied to
the complete, intermediate or the simple mixed meta-models (\ref{eq:completeMixedMetaModel}), (\ref{eq:intermediateMixedMetaModel}) and (\ref{eq:simpleMixedMetaModel}), respectively. This is given by 
\cite{Kuhn2005} that we briefly recall, without detailing their assumptions (M1)-(M5) and (SAEM1)-(SAEM4). 

\begin{proposition}[Kuhn and Lavielle]\label{prop:Kuhn}

Under assumptions (M1)-(M5) and (SAEM1)-(SAEM4) for the complete, intermediate or  simple mixed meta-model, if the sequence $(s_k)$ stays in a compact set, the SAEM  algorithm  produces a  sequence $(\hat\theta^{(k)})_{k\geq 1}$ which converges to the (local) maximum of the approximate likelihood $p_D(\by; \theta)$, $\bar p_D(\by; \theta)$ or $\tilde p_D(\by; \theta)$, respectively.
 \end{proposition}

 Now we study the impact of the meta-model approximations on the likelihoods.
 Our goal is to obtain a uniform control on the distance between the likelihood of the exact model $p(\by;\btheta)$
and the likelihoods of the three mixed meta-models $p_D(\by;\btheta)$, $\bar p_D(\by;\btheta)$ and $\tilde p_D(\by;\btheta)$ as a function of the quality of the meta-model. %Proofs are gathered   in Appendix. 
We start by the simple mixed meta-model.
\begin{proposition}\label{prop:errorSimple}
Let us consider the likelihoods $p(\by; \theta)$ (\ref{eq:Likelihood}) of the mixed model
(\ref{eq:MixedModel}) and $\tilde p_D(\by; \theta)$ (\ref{eq:simpleMetaLikelihood})  of the simple mixed meta-model 
(\ref{eq:simpleMixedMetaModel}) associated to a minimax design $D$. 
Assume that the support of the distribution of $\psi$ is compact. Assume that the functions $f$ and $m_D$ are uniformly 
bounded on the support of the distribution of  $\psi$. 
Assume that $f$ lies in the RKHS associated with the kernel $K$ satisfying to the same hypotheses as in Proposition \ref{prop:errorMetaModel}.
Then, there exists a constant $\tilde C_y$ which depends only on $\by$ such that
$$
 |p(\by; \theta)-\tilde p_D(\by; \theta)|\le  \tilde C_y\frac{n_{tot}}{\sigma_\varepsilon^{n_{tot}+2}}  G_K(a_D)\,
$$
where the function $G_K$ tends to $0$ when $a\rightarrow 0$ (defined in  Proposition \ref{prop:errorMetaModel}) and the constant $a_D$  is the covering distance of the design of experiments $D$. 
% il faut exiger dans les conditions que le noyau K est tel que G existe...

\end{proposition}
Recall that, when using a Gaussian kernel $K(x,x')$ for the meta-model approximation, the function $G_K$ is defined by $G_K(a)=Ce^{-\delta/a^2}$. Then, to ensure that this covering distance is small, we need a global upper-bound, uniformly in $\psi$. This is true when  the support of the distribution of $\psi$
is compact. Under this assumption, we obtain that the covering distance $G_K(a_D)$ can be as small as required provided there is a sufficient number of points $n_D$ in the design.
Thus providing a rich design $D$ during the pre-computation step allows to control as finely as we want the error induced on the likelihoods. 

 Now, we can study the distance between the three mixed meta-models. \begin{proposition}\label{prop:errorComplete}

Let us consider the likelihoods $p_D(\by; \theta)$ (\ref{eq:completeMetaLikelihood}) 
of the complete mixed meta-model (\ref{eq:completeMixedMetaModel}),  $\bar p_D(\by; \theta)$ 
(\ref{eq:intermediateMetaLikelihood})  of the intermediate mixed meta-model (\ref{eq:simpleMixedMetaModel})  and $\tilde p_D(\by; \theta)$ (\ref{eq:simpleMetaLikelihood})  of the simple mixed meta-model (\ref{eq:simpleMixedMetaModel}) associated to a minimax design $D$.  

Under the same hypotheses as Proposition \ref{prop:errorSimple}, there exist two constants $  C_y$ and $\bar C_y$which depend only on $\by$ such that
\begin{eqnarray*}
 |p_D(\by;\theta)-\tilde p_D(\by)|&\le&  C_y \frac{n_{tot}}{ \sigma_\varepsilon^{n_{tot}+2}}G_K(a_D), \\
 |\bar p_D(\by;\theta)-\tilde p_D(\by)|&\le&\bar C_y \frac{n_{tot}}{ \sigma_\varepsilon^{n_{tot}+2}}G_K(a_D).
 \end{eqnarray*}

\end{proposition}
Therefore, this  guarantees a control between the likelihood of any of the 
mixed meta-model and the likelihood of the exact mixed model. 

 With regularity hypotheses on the Hessian matrix  of each likelihood, results similar to \cite{Donnet2007}  can be obtained: The distance between the maximum of the exact likelihood $p(\by; \theta)$ and the maximum of the approximate likelihoods $p_D(\by; \theta)$, $\bar p_D(\by; \theta)$ or $\tilde p_D(\by; \theta)$
 can be as small as we want, as soon as the design $D$ is rich enough. 
%table avec les 25
\begin{table*}
\begin{tabular}{lrrrrrrrrrrr}
\hline
 Parameter && \multicolumn{3}{c}{Intermediate }& &\multicolumn{3}{c}{Simplified   }&& {Exact } \\
&  & \multicolumn{3}{c}{ meta-model} && \multicolumn{3}{c}{  meta-model} &&model\\
\hline
& $n_D$&25& 50 &100  &&25&50&100 && \\
\hline
\multirow{3}{*}{$\mu_{\log V}$}&Bias &-0.508 &-0.025   &  0.121  &&-0.483 &  -0.015&   0.089 && -0.390\\
                      &RMSE&  0.048&0.043   & 0.042 && 0.048& 0.042 & 0.041  &&  0.063  \\
                      &Cov.&92.8 &93.3  &  93.3  && 92.1 & 93.6 &   93.8   &&  86.3 \\
\hline                      
\multirow{3}{*}{$\mu_{\log k_a}$}&Bias &-4.920 & -1.389   &  -0.794  && -4.872 &  -1.396 &   -0.904&& -2.079\\
                      &RMSE & 0.860&0.545   &   0.476  && 0.870&  0.541 &  0.502 &&  0.797 \\
                      &Cov.  &79.2 &86.9   &  89.1  && 77.8 & 85.9 &    87.4 &&   81.0 \\                     
\hline
\multirow{3}{*}{$\mu_{\log V_m }$}&Bias& -2.067 & -0.599  & 0.014  &&-1.930 &  -0.577  &    -0.138  && -1.566\\
                      &RMSE& 0.392&  0.333   &  0.314 &&0.401 &  0.328&   0.327  &&  0.680\\
                      &Cov.& 87.3& 88.1   & 88.9  && 85.5& 88.9 &   89.2   &&   82.5  \\                      
\hline
\multirow{3}{*}{$\omega_{\log V }^2$}&Bias & 4.569& -2.408   &  -2.185      &&4.270 & -2.108  &  -2.215 &&  -4.393\\
                                 &RMSE & 6.487 &  5.526  & 5.276  && 6.359 & 5.445 &  5.287  && 9.461\\
                                 &Cov. &94.5  & 92.2  &  91.9 && 94.8& 93.3 &  93.5  &&  84.8 \\                      
\hline
\multirow{3}{*}{$\omega_{\log k_a}^2$}&Bias & 1.755&  -3.822   &   -6.797    &&1.935 &  -5.022  &    -6.799  &&  -1.705\\
                                 &RMSE &  17.398&  16.382   &  16.590  &&17.465  & 16.305 &  16.870  && 23.416 \\
                                 &Cov. & 84.4 & 82.6   &  81.4  && 84.9& 81.6 &   80.4&& 80.0 \\                      
\hline
\multirow{3}{*}{$\omega_{\log V_m}^2$}&Bias&  -33.721& -30.387  &   -30.408    && -33.916 & -29.946&   -30.148&&  1.867  \\
                                 &RMSE&  15.975 &13.236   & 13.039  && 15.981 & 12.914 &  13.236 &&  17.039\\
                                 &Cov. & 62.6 & 65.6  & 65.8  && 62.3 & 69.6 &    67.7 &&  83.8 \\                      
\hline
\multirow{2}{*}{$\sigma^2_{\epsilon}$} &Biais& 2.449&1.975   &   2.337    && 5.054 &  2.648 &   2.450 &&  0.308\\
                                       &RMSE&  0.354 &  0.302    &  0.370      &&  0.650 & 0.426&  0.397 &&  0.177\\
\hline
\end{tabular} 
\caption{Michaelis-Menten pharmacokinetic simulations: relative bias ($\%$), relative MSE ($\%$)  and coverage rate ($\%$)  computed over 1000   simulations,
with the intermediate meta-, the simple meta- and the exact mixed models.
Meta-models are built with either $n_D=25$, $n_D=50$  or $n_D=100$ design points. Coverage rate (Cov.) is the coverage rate
of the $95\%$ confidence interval based on the stochastic approximation of the Fisher matrix. \label{tabMM:varestimated}}
\end{table*}

%%%%%%%%%%%%% 
 
\section{Simulation study}\label{sec:simulation}
The objective of this study is to compare the main statistical properties of the estimation with  the   mixed meta-models and compare them to the exact mixed model. Two examples are illustrated  below, using standard ODE pharmacokinetics (PK) models.

 \begin{table*}
 %table avec a estime par fminsearch
\begin{tabular}{lrrrrrrrrr}
\hline
 Parameter && \multicolumn{2}{c}{Intermediate }& &\multicolumn{2}{c}{Simple   }&& {Exact } \\
&  & \multicolumn{2}{c}{ meta-model} && \multicolumn{2}{c}{  meta-model} &&model\\
\hline
& $n_D$& 50 &100  &&50&100 && \\
\hline
\multirow{3}{*}{$\mu_{\log k_e}$}&Bias   &0.101   & 0.007   &&   -0.320 &    0.007&&  0.003\\
                      &RMSE   & 0.004  &  0.005  &&   0.005 & 0.005 &&  0.005\\
                      &Cov.    & 94.2   & 94.4   &&   90.6 &    94.6 &&  93.9\\
\hline                      
\multirow{3}{*}{$\mu_{\log k_a}$}&Bias    & -2.441   & 0.001   &&      -8.380 &    0.008&&  -0.220\\
                      &RMSE  & 0.222   & 0.162  &&   0.910 &   0.160 &&  0.160\\
                      &Cov.  &  90.9  &  95.6  &&   59.6 &  95.3   &&  95.6 \\                     
\hline
\multirow{3}{*}{$\mu_{\log C_l}$}&Bias  & 0.388  & 0.036   &&     0.160   &    0.036  &&  -0.004\\
                      &RMSE   & 0.004  & 0.003  &&   0.003 &   0.003 &&  0.003\\
                      &Cov.  & 87.6   & 95.1  &&   93.4 &    94.7   &&  94.9\\                      
\hline
\multirow{3}{*}{$\omega_{\log k_e}^2$}&Bias   & -12.113   &  -2.745     &&  -23.200  &  -2.780 && -3.400\\
                                 &RMSE  & 7.131   & 6.404  &&   9.730 &   6.530  &&  6.460 \\
                                 &Cov.  & 83.2  & 91.5   &&   65.7 &  90.5    &&  90.3 \\                      
\hline
\multirow{3}{*}{$\omega_{\log k_a}^2$}&Bias  & -20.485   &  -3.442     &&  20.900  &    -3.320  &&  -2.440\\
                                 &RMSE  & 10.696  &  5.911  &&   13.500 &    5.930  &&  6.050\\
                                 &Cov.  &  72.3  &  89.7  &&   96.9 &    89.2   &&  90.2\\                      
\hline
\multirow{3}{*}{$\omega_{\log C_l}^2$}&Bias  & 0.375   &    -1.145   &&    -8.100 &   -1.100 &&  -2.660\\
                                 &RMSE   & 5.944   & 5.726  &&   5.810 &    5.690 &&  5.650\\
                                 &Cov.  &  92.6  &  92.0  &&   87.5 &    92.8  &&  91.1 \\                      
\hline
\multirow{2}{*}{$\sigma^2_{\epsilon}$} &Biais &  -45.262   &   -0.612     &&    16.000 &   -0.009 && -0.023\\
                                       &RMSE  & 20.719    &  0.232      &&    2.950 &   0.220 &&  0.220\\
\hline
\end{tabular}
%\caption{Simulation results when $\sigma^2_{\epsilon}$
%is estimated (THANKS TO the computation of a sufficient statistic for $\sigma^2_{\epsilon}$): relative bias and 
\caption{One compartment simulations: relative bias ($\%$), relative MSE ($\%$)  and coverage rate ($\%$)   computed over 1000   simulations,
with the intermediate meta-, the simple meta- and the exact mixed models.
Meta-models are built with either $n_D=50$  or $n_D=100$ design points. Coverage rate (Cov.) is the coverage rate
of the $95\%$ confidence interval based on the stochastic approximation of the Fisher matrix. \label{tab:varestimated}}
\end{table*}

\subsection{Michaelis-Menten pharmacokinetic model}

\subsubsection{Simulation settings}\label{pkMM}
Let us now consider a one-compartment pharmacokinetic model, first order absorption and Michaelis Menten elimination. A dose $D$ of a drug is given to a patient by intra-venous bolus. 
The concentration of the drug in the body along time in then described by the following ordinary differential equation
$$ \frac{df}{dt} = - \frac{{V_m} \cdot f}{k_m +f} +k_a\cdot \frac{D}{V}\cdot \exp(-k_a\cdot t), \quad f(t_0) = 0$$ where $V$ is the volume of distribution, $V_m$ is the maximum elimination rate (in amount per time unit), $k_m$ is the Michaelis-Menten constant (in concentration unit) and
$k_a$ is the absorption constant.  
We consider $\log k_m$ as fixed to $-2.5$. The individual parameter $\psi$ consists in $\log V$, $\log k_a$ and $\log V_m$.
We assume a Gaussian distribution on the logarithm of these parameters with mean $(\mu_{\log V},\mu_{\log k_a},\mu_{\log V_m})=(2.5,1,-0.994)$
and a diagonal covariance matrix with terms $(\omega^2_{\log V},\omega^2_{\log k_a},\omega^2_{\log V_m})=(0.09,0.09,0.09)$.
Then a homoscedastic additive error model is simulated with a standard error $\sigma_\varepsilon = 0.1$. 

We implement the four algorithms: SAEM on the original mixed model, SAEM on the complete, intermediate and simple mixed meta-model. 
 For the   meta-model SAEM algorithms, we use successively   $n_D=25$, $n_D=50$ and $n_D=100$   number of points in the design of experiments for the Gaussian process emulator. 
The covariance is Gaussian, and the regression functions $H$ are linear functions. More sophisticated choices in the regression functions and in the
kernel can be made. However, our goal in this section is to illustrate in a quite simple case the efficiency of 
the combination of the Gaussian process emulation with the SAEM-MCMC algorithm.
 In the pre-computation step,  for a given value $\psi$,  the ODE solver provides $f(t, \psi)$ for each time of measurement. 
Thus, the design of numerical experiment is only built over the values of $\psi$ and not $t$. We have compared two approaches one where $t$ is
considered as an additional input and the other where a meta-model is built for each time $t$. Based on the comparison of the
quality of the approximations, we have kept the second one which is quite simple to deal with. However, more sophisticated approaches
can be tested \citep[see][for a review]{rougier:2008}. %on n'en parle que ici suffisant ou le dire aussi dans la partie spécifique au métamodèle ?
%The same SAEM algorithms were run as in the first example. 
%The Gaussian process emulators were also built with $n_D=25$, $n_D=50$ and $n_D=100$ points, linear
%regression functions and a Gaussian covariance kernel, in the same fashion as in subsection \ref{numtheo}.
The approximation is built over the domain: $[1.6 ;3.3]$ for $\log V$, $[0; 2.1]$ for $\log k_a$ and $[-1.6; -0.3]$ for $\log V_m$.
The starting values for the parameters are 
$\widehat\mu_{\log V}^{(0)} = 2$, $\widehat\mu_{\log k_a}^{(0)}=0.5$, $\widehat\mu_{\log V_m}^{(0)}=-0.5$, $\widehat\omega_{\log V}^{2(0)}=\widehat\omega_{\log k_a}^{2(0)}=\widehat\omega_{\log V_m}^{2(0)}=0.1$  and $\widehat\sigma_\varepsilon^{(0)}=0.3$.

\subsubsection{Results}
The computation times for one run of SAEM (100 iterations of SAEM, with 15 iterations of MCMC at each SAEM iteration) were the following:
around 15 min for the exact mixed model (requiring solving the ODE at each iteration of MCMC), 
around 30 min for the complete mixed meta-model with $n_D=50$ (requiring inverting the $\bC(\bt,\bpsi)$ at each iteration of MCMC),
around 80 sec for   the intermediate  model and 30 sec for the simple one.
Therefore, in the following, we only present the results for the exact mixed model (as a benchmark) and the intermediate and simple mixed meta-models.

Relative bias and relative root mean square error (RMSE) are computed for each population parameter  from 1000 replications and presented in Table \ref{tabMM:varestimated}. 
The $95\%$ coverage rates correspond to the coverage rate 
of the  confidence interval on parameters based on the stochastic approximation of the Fisher Information matrix. 
In this example, the two meta-models have   good performances even with only $25$ points in the design of experiments and increasing $n_D$ decreases the bias.  The parameter $\omega_{\log V_m}$ is biased when using a meta-model (whatever $n_D$) while it is not with the exact model. This may be due to the error of approximation that is not completely taken into account.

\subsection{First order pharmacokinetic model}
\label{numtheo}
\subsubsection{Simulation settings}
Let us consider a one-compartment PK model with first order absorption and elimination, with a dose $D$ of a drug. The concentration of the drug in the body along time in then described by the following ordinary differential equation
$$ \frac{df}{dt} = D \frac{k_ak_e}{C_l}exp(-k_at) -k_ef, \quad f(t_0) = 0$$
where $k_a$ and $k_e$ are the absorption and elimination constants, $C_l$ is the clearance. 
We consider the PK parameters of theophyllin  \citep{Pinheiro2000}:
$\log k_e = -2.52$, $\log k_a=0.4$, $\log C_l=-3.22$. One dataset of 36 patients is simulated with a dose D=6 $mmol$ and measurements at time $t =0.25, 0.5, 1, 2, 3.5, 5, 7, 9, 12$ hours.
The random effects 
were simulated assuming a Gaussian distribution for the logarithm of the parameters with a diagonal variance-covariance matrix $\Omega$ with the following diagonal elements: $\omega_{\log k_e }^2=\omega_{\log k_a }^2=\omega_{\log C_l}^2=0.01$.  
Then a homoscedastic additive error model is simulated with a standard error $\sigma_\varepsilon = 0.1$.

The same SAEM algorithms were run as in the first example. The Gaussian process emulators were built with $n_D=50$ and $n_D=100$ points, linear
regression functions and a Gaussian covariance kernel, in the same fashion as in subsection \ref{pkMM}.
The domain where the approximation is built is $[-4; -1]$ for $\log k_e$, $[0; 2]$ for $\log k_a$ and $[-4.5; 2]$ for $\log C_l$.
% The covariance is Gaussian, and the regression functions $H$ are linear functions. More sophisticated choices in the regression functions and in the
% kernel can be made. However, our goal in this section is to illustrate in a quite simple case the efficiency of 
% the combination of the Gaussian process emulation with the SAEM-MCMC algorithm.
%DONNER LES VALEURS DE DEPART DE L'ALGO ???
The starting values for the parameters are 
$\mu_{\log k_e}^{(0)} = -3$, $\mu_{\log k_a}^{(0)}=1$, $\mu_{\log C_k}^{(0)}=-3$, $\omega_{\log k_e}^{2(0)}=\omega_{\log k_a}^{2(0)}=\omega_{\log C_l}^{2(0)}=0.1$  and $\sigma_\varepsilon^{(0)}=0.3$.

\subsubsection{Results}
The computation times are the same as before. 
Relative bias, relative RMSE and coverage rate   computed from 1000 replications are presented in Table \ref{tab:varestimated}. 
When the design of numerical experiments in the pre-computation step
has $100$ points, the estimates obtained with the mixed meta-models have similar performance
to the ones obtained with the exact mixed model.
With only $50$ points in the design, 
the estimates with the mixed meta-models are less accurate, especially $\sigma_\varepsilon$ with the intermediate
mixed model. 
There is a clear improvement of the quality of the estimates with the intermediate mixed meta-model for the parameters concerning the means.
Recall that the simple mixed meta-model neglects the approximation of the function $f$. Therefore, 
 taking into account the errors of the approximation of the Gaussian process emulator  in the model
prevents from a systematic
bias in the estimates. 
However, since the correlation between the Gaussian process emulator approximation errors are set to zero for the sake of simplicity, the estimation of   $\sigma_\varepsilon$ may be less accurate but usually this parameter if of less interest.

\section{Concluding remarks}\label{sec:conclusion}

In the case of a mixed model where the regression function is a non-analytical solution of an ODE or of a PDE, we proposed to build a so-called
meta-model which is obtained thanks to a pre-computation step. It consists in running the ODE/PDE solver on a well chosen design of numerical experiments.
Once this meta-model is obtained, we use it as a surrogate of the regression function in the estimation procedure which is based on a SAEM-MCMC algorithm.
We derived three mixed meta-models depending on whether the additional source of uncertainty due to the approximation by the meta-model is taken
into account totally, partially or not at all (complete, intermediate and simple mixed meta-model).
In the complete mixed meta-model, there is a full covariance matrix accounting for dependencies induced by the meta-modeling errors
which slows down the SAEM-MCMC algorithm. That is why we have renounced to test it in the simulation study.
Further works are needed to design MCMC algorithm adapted to this case of non-independent individuals.
In the intermediate and simple mixed meta-model, the individuals are still independent
thus the SAEM-MCMC algorithm does not suffer from any computational burden.
We showed examples where even with a very few design points for the meta-model approximation, the estimation results are very satisfactory. 
We also showed an example where the intermediate meta-model improves the quality of the estimates of the parameters especially those accounting for 
the mean of the population parameters. 
\\

Since the quality of the approximation provided by the meta-model directly depends on the density of the numerical design of experiments
in the neighborhood of the input where the approximation is made, a sequential strategy for building an adaptive design reinforcing the meta-model
where the SAEM-MCMC algorithm identifies likely region for the parameters should improve the estimates. However,
this strategy would make the Markov property in the SAEM-MCMC procedure no longer to be valid. Therefore, there are theoretical 
questions which will be interesting to solve in order to ensure guarantees in this case.

 \section*{Acknowledgements}
Adeline Samson has been supported by   the LabEx \\
PERSYVAL-Lab (ANR-11-LABX-0025-01).   Les recherches menant aux présents résultats ont bénéficié d'un soutien financier du septième programme-cadre de l'Union européenne (7ePC/2007-2013) en vertu de la convention de subvention n 266638.
%\noindent The authors thank Marc Lavielle for his advices, suggestions and support.

\bibliographystyle{apalike}      % basic style, author-year citations
\bibliography{biblio2}   % name your BibTeX data base

% \begin{thebibliography}{}
% %
% % and use \bibitem to create references. Consult the Instructions
% % for authors for reference list style.
% % Format for Journal Reference
% Author, Article title, Journal, Volume, page numbers (year)
% Author, Book title, page numbers. Publisher, place (year)
% \end{thebibliography}

\section*{Appendix}
\subsection*{Proof of proposition \ref{prop:errorSimple}}
We have 
$$|p(\by;\theta)-\tilde p_D(\by;\theta)|\le \int |p(\by|\bpsi;\btheta)-\tilde p_D(\by|\bpsi;\btheta)|p(\bpsi)d\bpsi\,.$$
Therefore, we start by studying $|p(\by|\bpsi;\btheta)-\tilde p_D(\by|\bpsi;\btheta)|$: 
 
%\begin{tiny}
{\small
\begin{eqnarray*}
 &&(2\pi \sigma_\varepsilon^2)^{n_{tot}/2} |p(\by|\bpsi;\btheta)-\tilde p_D(\by|\bpsi;\btheta)|\\
 &=& \Bigg|
 \exp\bigg(-\frac{1}{2\sigma_\varepsilon^2}\sum_{ij}(y_{ij}-f(t_{ij},\psi_i))^2\bigg)
 -\exp\bigg(-\frac{1}{2\sigma_\varepsilon^2}\sum_{ij}(y_{ij}-m_D(t_{ij},\psi_i))^2\bigg)\Bigg|\\
&=&  \exp\left(-\frac{1}{2\sigma_\varepsilon^2}\sum_{ij}(y_{ij}-f(t_{ij},\psi_i))^2\right)\\
&&\times
\Bigg|1-\exp\bigg(-\frac{1}{2\sigma_\varepsilon^2}\sum_{ij}\Big((y_{ij}-m_D(t_{ij},\psi_i))^2
-(y_{ij}-f(t_{ij},\psi_i))^2\Big)\bigg) \Bigg|\\
&\le&  \Bigg|1-\exp\bigg(-\frac{1}{2\sigma_\varepsilon^2}\sum_{ij}\Big(f(t_{ij},\psi_i)-m_D(t_{ij},\psi_i))
(2y_{ij}-f(t_{ij},\psi_i)
-m_D(t_{ij},\psi_i)\Big)\bigg)\Bigg|\,.
 \end{eqnarray*}
} 
%\end{tiny}
Under the assumption that the functions $f$ and $m_D$ are uniformly bounded on the support of $\psi$,  there exists a constant $C_y$  which is uniform according to 
$\psi$ such that   $|2y_{ij}-f(t,\psi)-m_D(t,\psi)|\le C_y$. Proposition \ref{prop:errorMetaModel} implies that the approximation error due to the metamodel $|f(t_{ij},\psi_i)-m_D(t_{ij},\psi_i)|$ is controlled by inequality (\ref{eq:bornedeterministe}):
$$|f(t_{ij},\psi_i)-m_D(t_{ij},\psi_i)|\le \Vert f\Vert_{\mathcal{H}_K}G_K(a_D)\,.$$
Then there exists a constant $C_y$ depending only on $\by$ such that 
$$
(2\pi \sigma_\varepsilon^2)^{n_{tot}/2} |p(\by|\bpsi;\btheta)-\tilde p_D(\by|\bpsi;\btheta)|\leq  C_y \frac{n_{tot}}{2\sigma_\varepsilon^2}  \Vert f\Vert_{\mathcal{H}_K}G_K(a_D).
$$
Finally 
$$ |p(\by;\theta)-\tilde p_D(\by;\theta)|\le\frac{C_y}{(2\pi \sigma_\varepsilon^2)^{n_{tot}/2}}\frac{n_{tot}}{2\sigma_\varepsilon^2}  \Vert f\Vert_{\mathcal{H}_K}G_K(a_D)
\,.$$
$\Box$

\subsection{Proof of proposition \ref{prop:errorComplete}}
We study the distance between the two likelihoods $p_D$ and $\tilde p_D$. As in Proposition \ref{prop:errorSimple}, we start by studying $|p_D(\by|\bpsi;\btheta)-\tilde p_D(\by|\bpsi;\btheta)|$. We consider two Gaussian distributions with same expectations and different covariance matrix. Thus this distance is maximum when $\sum (y_{ij}-m_D(t_{ij},\psi_i))^2=0$. This yields
\begin{eqnarray*}
 &&(2\pi)^{n_{tot}/2} |p_D(\by|\bpsi;\btheta)-\tilde p_D(\by|\bpsi;\btheta)|\\
 &\leq& \left| \frac{1}{\sigma_\varepsilon^{n_{tot}}}-\frac{1}{\sqrt{|\sigma_\varepsilon^2I_{n_{tot}}+\bC_D|}}\right|\\
 &=&\frac{1}{ \sigma_\varepsilon^{n_{tot}}}\left|1-\frac{\sigma_\varepsilon^{n_{tot}}}{|\sigma_\varepsilon^2I_{n_{tot}}+\bC_D|^{1/2}}\right|\\
 &\le&\frac{1}{ \sigma_\varepsilon^{n_{tot}}}\left|1-\frac{\sigma_\varepsilon^{n_{tot}}}{(\sigma_\varepsilon^2+\frac{1}{n_{tot}}\sum_{ij} C_D(t_{ij},\psi_i;t_{ij},\psi_i))^{n_{tot}/2}}\right|\\
  &\le&\frac{1}{ \sigma_\varepsilon^{n_{tot}}}\left|1-\frac{1}{(1+\frac{1}{\sigma_\varepsilon^2n_{tot}}\sum_{ij} C_D(t_{ij},\psi_i; t_{ij},\psi_i))^{n_{tot}/2}}\right|\\
 \end{eqnarray*}
 where we use that the determinant, as a product of eigen values,  is smaller than a function of the trace of the matrix. Thus,  the sum is over the diagonal of the matrix $\bC_D$ i.e. the sum of the variances.
 Then, we obtain that there exists a constant $C$ such that 
 \begin{eqnarray*}
   |p_D(\by|\bpsi;\btheta)-\tilde p_D(\by|\bpsi;\btheta)|
   &\le& C \frac{1}{ \sigma_\varepsilon^{n_{tot}}}\left| \frac1{2\sigma_\varepsilon^2}\sum_{ij} C_D(t_{ij},\psi_i; t_{ij},\psi_i)\right|\\
   &\leq & C (2\pi)^{n_{tot}/2}\frac{n_{tot}}{ \sigma_\varepsilon^{n_{tot}+2}}G_K(a_D)
 \end{eqnarray*}
where the last inequality holds using Proposition \ref{prop:errorMetaModel}. 
Finally, we obtain
$$ |p_D(\by;\theta)-\tilde p_D(\by;\theta)|\le C_y \frac{n_{tot}}{ \sigma_\varepsilon^{n_{tot}+2}}G_K(a_D).$$
The proof is similar for the distance between the two likelihoods $\bar p_D$ and $\tilde p_D$.
$\Box$

\end{document}